\documentclass[12pt]{article}
\usepackage{latexsym}
\usepackage{amssymb}
\usepackage{amsmath}

\begin{document}


\newtheorem{theorem}{Theorem}[section]
\newtheorem{problem}{Problem}
\newtheorem{definition}{Definition}[section]
\newtheorem{lemma}{Lemma}[section]
\newtheorem{proposition}{Proposition}[section]
\newtheorem{corollary}{Corollary}[section]
\newtheorem{example}{Example}
\newtheorem{conjecture}{Conjecture}
\newtheorem{algorithm}{Algorithm}
\newtheorem{exercise}{Exercise}
\newtheorem{remarkk}{Remark}[section]

\newcommand{\be}{\begin{equation}}
\newcommand{\ee}{\end{equation}}
\newcommand{\bea}{\begin{eqnarray}}
\newcommand{\eea}{\end{eqnarray}}
\newcommand{\beq}[1]{\begin{equation}\label{#1}}
\newcommand{\eeq}{\end{equation}}
\newcommand{\beqn}[1]{\begin{eqnarray}\label{#1}}
\newcommand{\eeqn}{\end{eqnarray}}
\newcommand{\beaa}{\begin{eqnarray*}}
\newcommand{\eeaa}{\end{eqnarray*}}
\newcommand{\req}[1]{(\ref{#1})}

\newcommand{\lip}{\langle}
\newcommand{\rip}{\rangle}
\newcommand{\uu}{\underline}
\newcommand{\oo}{\overline}
\newcommand{\La}{\Lambda}
\newcommand{\la}{\lambda}
\newcommand{\eps}{\varepsilon}
\newcommand{\om}{\omega}
\newcommand{\Om}{\Omega}
\newcommand{\ga}{\gamma}
\newcommand{\rrr}{{\Bigr)}}
\newcommand{\qqq}{{\Bigl\|}}

\newcommand{\dint}{\displaystyle\int}
\newcommand{\dsum}{\displaystyle\sum}
\newcommand{\dfr}{\displaystyle\frac}
\newcommand{\bige}{\mbox{\Large\it e}}
\newcommand{\integers}{{\Bbb Z}}
\newcommand{\rationals}{{\Bbb Q}}
\newcommand{\reals}{{\rm I\!R}}
\newcommand{\realsd}{\reals^d}
\newcommand{\realsn}{\reals^n}
\newcommand{\NN}{{\rm I\!N}}
\newcommand{\DD}{{\rm I\!D}}
\newcommand{\LL}{{\rm I\!L}}
\newcommand{\degree}{{\scriptscriptstyle \circ }}
\newcommand{\dfn}{\stackrel{\triangle}{=}}
\def\complex{\mathop{\raise .45ex\hbox{${\bf\scriptstyle{|}}$}
     \kern -0.40em {\rm \textstyle{C}}}\nolimits}
\def\hilbert{\mathop{\raise .21ex\hbox{$\bigcirc$}}\kern -1.005em {\rm\textstyle{H}}} 
\newcommand{\RAISE}{{\:\raisebox{.6ex}{$\scriptstyle{>}$}\raisebox{-.3ex}
           {$\scriptstyle{\!\!\!\!\!<}\:$}}} 

\newcommand{\hh}{{\:\raisebox{1.8ex}{$\scriptstyle{\degree}$}\raisebox{.0ex}
           {$\textstyle{\!\!\!\! H}$}}}

\newcommand{\OO}{\won}
\newcommand{\calA}{{\cal A}}
\newcommand{\calB}{{\cal B}}
\newcommand{\calC}{{\cal C}}
\newcommand{\calD}{{\cal D}}
\newcommand{\calE}{{\cal E}}
\newcommand{\calF}{{\cal F}}
\newcommand{\calG}{{\cal G}}
\newcommand{\calH}{{\cal H}}
\newcommand{\calK}{{\cal K}}
\newcommand{\calL}{{\cal L}}
\newcommand{\calM}{{\cal M}}
\newcommand{\calO}{{\cal O}}
\newcommand{\calP}{{\cal P}}
\newcommand{\calR}{{\cal R}}
\newcommand{\calX}{{\cal X}}
\newcommand{\calXX}{{\cal X\mbox{\raisebox{.3ex}{$\!\!\!\!\!-$}}}}
\newcommand{\calXXX}{{\cal X\!\!\!\!\!-}}
\newcommand{\gi}{{\raisebox{.0ex}{$\scriptscriptstyle{\cal X}$}
\raisebox{.1ex} {$\scriptstyle{\!\!\!\!-}\:$}}}
\newcommand{\intsim}{\int_0^1\!\!\!\!\!\!\!\!\!\sim}
\newcommand{\intsimt}{\int_0^t\!\!\!\!\!\!\!\!\!\sim}
\newcommand{\pp}{{\partial}}
\newcommand{\al}{{\alpha}}
\newcommand{\sB}{{\cal B}}
\newcommand{\sL}{{\cal L}}
\newcommand{\sF}{{\cal F}}
\newcommand{\sE}{{\cal E}}
\newcommand{\sX}{{\cal X}}
\newcommand{\R}{{\rm I\!R}}
\newcommand{\vp}{\varphi}
\newcommand{\N}{{\rm I\!N}}
\def\ooo{\lip}
\def\ccc{\rip}
\newcommand{\ot}{\hat\otimes}
\newcommand{\rP}{{\Bbb P}}
\newcommand{\bfcdot}{{\mbox{\boldmath$\cdot$}}}

\renewcommand{\varrho}{{\ell}}
\newcommand{\dett}{{\textstyle{\det_2}}}
\newcommand{\sign}{{\mbox{\rm sign}}}
\newcommand{\TE}{{\rm TE}}
\newcommand{\TA}{{\rm TA}}
\newcommand{\E}{{\rm E\,}}
\newcommand{\won}{{\mbox{\bf 1}}}
\newcommand{\Lebn}{{\rm Leb}_n}
\newcommand{\Prob}{{\rm Prob\,}}
\newcommand{\sinc}{{\rm sinc\,}}
\newcommand{\ctg}{{\rm ctg\,}}
\newcommand{\loc}{{\rm loc}}
\newcommand{\trace}{{\,\,\rm trace\,\,}}
\newcommand{\Dom}{{\rm Dom}}
\newcommand{\ifff}{\mbox{\ if and only if\ }}
\newcommand{\proof}{\noindent {\bf Proof:\ }}
\newcommand{\remark}{\noindent {\bf Remark:\ }}
\newcommand{\remarks}{\noindent {\bf Remarks:\ }}
\newcommand{\note}{\noindent {\bf Note:\ }}

\newcommand{\boldx}{{\bf x}}
\newcommand{\boldX}{{\bf X}}
\newcommand{\boldy}{{\bf y}}
\newcommand{\boldR}{{\bf R}}
\newcommand{\uux}{\uu{x}}
\newcommand{\uuY}{\uu{Y}}

\newcommand{\limn}{\lim_{n \rightarrow \infty}}
\newcommand{\limN}{\lim_{N \rightarrow \infty}}
\newcommand{\limr}{\lim_{r \rightarrow \infty}}
\newcommand{\limd}{\lim_{\delta \rightarrow \infty}}
\newcommand{\limM}{\lim_{M \rightarrow \infty}}
\newcommand{\limsupn}{\limsup_{n \rightarrow \infty}}

\newcommand{\ra}{ \rightarrow }

\newcommand{\ARROW}[1]
  {\begin{array}[t]{c}  \longrightarrow \\[-0.2cm] \textstyle{#1} \end{array} }

\newcommand{\AR}
 {\begin{array}[t]{c}
  \longrightarrow \\[-0.3cm]
  \scriptstyle {n\rightarrow \infty}
  \end{array}}

\newcommand{\pile}[2]
  {\left( \begin{array}{c}  {#1}\\[-0.2cm] {#2} \end{array} \right) }

\newcommand{\floor}[1]{\left\lfloor #1 \right\rfloor}

\newcommand{\mmbox}[1]{\mbox{\scriptsize{#1}}}

\newcommand{\ffrac}[2]
  {\left( \frac{#1}{#2} \right)}

\newcommand{\one}{\frac{1}{n}\:}
\newcommand{\half}{\frac{1}{2}\:}

\def\le{\leq}
\def\ge{\geq}
\def\lt{<}
\def\gt{>}

\def\squarebox#1{\hbox to #1{\hfill\vbox to #1{\vfill}}}
\newcommand{\qed}{\hspace*{\fill}
           \vbox{\hrule\hbox{\vrule\squarebox{.667em}\vrule}\hrule}\bigskip}

\title{Solution of the Monge-Amp\`ere Equation on Wiener Space for
  log-concave measures}

\author{D. Feyel and A. S. \"Ust\"unel}
\date{ }

\maketitle
\begin{abstract}
\noindent
In this work we prove that the unique  $1$-convex solution of the
Monge-Kantorovitch measure transportation problem between the Wiener
measure and a target measure which has a log-concave density w.r.to
the Wiener measure is also the strong solution of   the Monge-Amp\`ere
equation in the frame of infinite dimensional Fr\'echet spaces. We
enhance also the polar factorization results  of the mappings which
transform a spread measure to another one in terms of the measure
transportation of Monge-Kantorovitch.
\end{abstract}

\noindent

\section{Introduction}

In 1781, G. Monge has launched his famous problem \cite{Monge}, which can be
expressed in terms of the modern mathematics as follows: given two
probability measures $\rho$ and $\nu$ on $\R^n$, find  the  map
$T:\R^n\to \R^n$ such that $T\rho=\nu$ \footnote{$T\rho$ means the
  image of the measure $\rho$ under the map $T$} and $T$ is also the solution
of the minimization problem 
\begin{equation}
\label{monge-prob}
\inf_U\left\{\int_{\R^n}c(x,U(x)) \rho(dx)\right\}\,,
\end{equation}
where the infimum is taken between all the maps $U:\R^n\to \R^n$ such
that $U\rho=\nu$ and where  $c:\R^n\times \R^n\to \R_+$ is a positive,
measurable function, called usually the cost function. In the original
problem of Monge, the cost function $c(x,y)$ was $|x-y|$ and the
dimension $n$ was three. Later other
costs have been considered, between them, the most popular one which
is also abundantly studied,  is the case where $c(x,y)=|x-y|^2$. After
several tentatives (cf., \cite{App-1,App-2}), in the 1940's this
highly nonlinear  problem of 
Monge has been reduced to a linear problem by Kantorovitch,
cf.\cite{Kan},  in the following way: let $\Sigma(\rho,\nu)$ be the set
of probability measures on $\R^n\times\R^n$, whose first marginals are
$\rho$ and the second marginals are $\nu$. Find the element(s) of
$\Sigma(\rho,\nu)$ which are the solutions of the minimization
problem:
\begin{equation}
\label{mkp-prob}
\inf_{\beta\in\Sigma(\rho,\nu)}\left\{\int_{\R^n\times\R^n}c(x,y)
d\beta(x,y)\right\}\,.
\end{equation}
It is obviuous that
$\Sigma(\rho,\nu)$ is a convex, compact set under the weak*-topology
of measures, hence, in case, the cost function $c$ has some regularity
properties, like being lower semi-continuous,  this problem would
have solutions. If any one of them is supported by the graph of a map
$T:\R^n\to \R^n$, then obviously, $T$ will be also a solution of the
original problem of Monge \ref{monge-prob}. Since that time, the
problem (\ref{mkp-prob}) is called the Monge-Kantorovitch problem (MKP).
The program of Kantorovitch  has been followed by several
people and a major contribution has been done  by Sudakov
\cite{Sud}. In the early 90's there has been another impetus to this
problem, cf., \cite{BRE}, where it has been discovered  the important
role played by the convex functions in the construction of  the
solutions  of the  MKP and of the problem of Monge (cf.,
\cite{Mc1,Mc2}). We refer the reader to   \cite{Feyel} and to
\cite{Vil} for  recent surveys.  

In \cite{F-U3}, we have solved
the MKP and the problem of Monge in the infinite dimensional case,
where  the measures are concentrated in  a Fr\'echet  space $W$
into which a Hilbert space $H$  is injected densely and continuously. 
We call $H$   the Cameron-Martin space in reference to the
Gaussian case. The cost function  is defined on $W\times W$ as
\beaa
c(x,y)&=&|x-y|_H^2 \mbox{ \rm{if} }x-y\in H\\
&=&\infty \mbox{ \rm{if} }x-y\notin H\,,
\eeaa
 where $|\cdot|_H$ denotes the Euclidean norm of
$H$. Because of this choice, in comparison to the finite dimensional 
space, the situation becomes quite singular, since, in general, the
Cameron-Martin space $H$ is a negligeable set (i.e., of null measure)
with respect to almost all  reasonable 
measures for which one can expect to have solutions of  the
 problems of Monge and of  MKP. On the other hand, due to
the potential applications to several problems of stochastic analysis
and physics, this cost function is particularly important. For
example, it is particularly well-adapted  to the  study of the
absolute continuity of the image of the Wiener measure under the
perturbations of identity, which is a subject under investigation
since the early works of N. Wiener, R.H. Cameron and W.T. Martin and of
several other mathematicians and engineers who have made worthy
contributions  (cf. the list of references of \cite{BOOK}). 


This paper is devoted to the applications and some furthere
developments of the subject. At first  we give a  generalization of
the polar factorization of  vector fields which map a  probability
measure on $W$ to another one such that one of them is spread   (cf. the
preliminaries) and  the two measures  are at finite Wasserstein distance
from each  other (without any absolute continuity hypothesis).  As an
example we treat in detail the case of the infinite dimensional
Gaussian measures.  
   
 The proof of the fact that the transport map, when the target measure
 has an $H$-log-concave density, satisfies  the
 functional analytic (or strong) Monge-Amp\`ere
 equation is probably the most important contribution of this
 paper. In \cite{F-U3}, we have
 studied the Monge-Amp\`ere equation for the upper and lower bounded
 densities with respect to the Wiener measure. The main difficulty in
 this infinite dimensional case stems from the lack of
 regularity of the transport 
 potentials, in fact we only know that these functions are in the
 Sobolev space $\DD_{2,1}$, i.e., they have only first order Sobolev
 derivatives. However, to write the Gaussian Jacobian, we need them to
 be  second order Sobolev derivatives taking values in the space of
 Hilbert-Schmidt operators on the Cameron-Martin space $H$. This
 difficulty is worse than those we encounter in the finite 
dimensional case, since in the latter the Hilbert-Schmidt property holds
automatically. Moreover,  in the  finite dimensional situation the
lack of second order derivatives is  solved  with the help of the 
Alexandroff derivatives  of the convex functions. In the infinite
dimensional case the situation is worse: the transport
potentials are not in general convex, nor $H$-convex (which is a
more reasonable  requirement than being convex, cf. \cite{F-U1}), but
only $1$-convex in the Cameron-Martin space direction. Hence their
second order derivatives in the sense of 
distributions are not in general measures; even if this happens in
some exceptional situations, their absolutely continuous parts do not
take values in the space of Hilbert-Schmidt operators, a condition
which is indispensable to write down the Jacobian of the transport
map. Hence it is impossible in general to construct the strong solutions of the
Monge-Amp\`ere equation.  
In Section \ref{st-MA}, combining  the finite dimensional
results of Caffarelli \cite{Caf} with Wiener space analysis, we solve
completely this problem when the target measure is
$H$-log-concave. More precisely, we show that 
the transport  potential has a second order derivative as an operator
valued map and then using some celebrated identity of Wiener space
analysis, we also prove that this second derivative  takes its values
in the space of Hilbert-Schmidt operators, hence we can write the corresponding
Jacobian which includes  the modified Carleman-Fredholm determinant,
cf. \cite{BOOK} and finally we prove that the transport potential is the
unique $1$-convex strong  solution of the Monge-Amp\`ere equation. 
 In Section \ref{Ito-section}  we show
that all these difficulties disappear if we use the natural Ito
Calculus and we can calculate  the It\^o  Jacobian (cf. Theorem \ref{M-A-Ito})
using the natural Brownian motion which is associated to the solution
of the Monge problem. In fact, with It\^o parametrization, the
complications are absorbed by the filtrations  of forward and backward
transport processes (i.e., maps). We give also the delicate relations
between the polar factorization of the absolutely continuous
transformations of the Wiener measure and the Brownian motions which
appear in the semimartingale decomposition of the transport process
with respect to its natural filtration. 

\section{Preliminaries and notations}
\label{preliminaries}
Let $W$ be a separable Fr\'echet space equipped with  a Gaussian
measure $\mu$ of zero mean whose support is the whole
space\footnote{The reader may assume that $W=C(\R_+,\R^d)$, $d\geq 1$
  or $W=\R^{\N}$.}. The
corresponding Cameron-Martin space is denoted by $H$. Recall that the
injection $H\hookrightarrow W$ is compact and its adjoint is the
natural injection $W^\star\hookrightarrow H^\star\subset
L^2(\mu)$. The triple $(W,\mu,H)$ is called 
an abstract Wiener space. Recall that $W=H$ if and only if $W$ is
finite dimensional. A subspace $F$ of $H$ is called regular if the
corresponding orthogonal projection 
has a continuous extension to $W$, denoted again  by the same letter.
It is well-known that there exists an increasing sequence of regular
subspaces $(F_n,n\geq 1)$, called total,  such that $\cup_nF_n$ is
dense in $H$ and in $W$. Let $V_n$  be the
$\sigma$-algebra generated by $\pi_{F_n}$, then  for any  $f\in
L^p(\mu)$, the martingale  sequence 
$(E[f|V_n],n\geq 1)$ converges to $f$ (strongly if 
$p<\infty$) in $L^p(\mu)$. Observe that the function
$f_n=E[f|V_n]$ can be identified with a function on the
finite dimensional abstract Wiener space $(F_n,\mu_n,F_n)$, where
$\mu_n=\pi_n\mu$. 

Since the translations of $\mu$ with the elements of $H$ induce measures
equivalent to $\mu$, the G\^ateaux  derivative in $H$ direction of the
random variables is a closable operator on $L^p(\mu)$-spaces and  this
closure will be denoted by $\nabla$ cf.,  for example
\cite{ASU}. The corresponding Sobolev spaces 
(the equivalence classes) of the  real random variables 
will be denoted as $\DD_{p,k}$, where $k\in \NN$ is the order of
differentiability and $p>1$ is the order of integrability. If the
random variables are with values in some separable Hilbert space, say
$\Phi$, then we shall define similarly the corresponding Sobolev
spaces and they are denoted as $\DD_{p,k}(\Phi)$, $p>1,\,k\in
\NN$. Since $\nabla:\DD_{p,k}\to\DD_{p,k-1}(H)$ is a continuous and
linear operator its adjoint is a well-defined operator which we
represent by $\delta$. In the case of classical Wiener space, i.e.,
when $W=C(\reals_+,\reals^d)$, then $\delta$ coincides with the It\^o
integral of the Lebesgue density of the adapted elements of
$\DD_{p,k}(H)$ (cf.\cite{ASU}). 

For any $t\geq 0$ and measurable $f:W\to \reals_+$, we note by
$$
P_tf(x)=\int_Wf\left(e^{-t}x+\sqrt{1-e^{-2t}}y\right)\mu(dy)\,,
$$
it is well-known that $(P_t,t\in \reals_+)$ is a hypercontractive
semigroup on $L^p(\mu),p>1$,  which is called the Ornstein-Uhlenbeck
semigroup (cf.\cite{ASU}). Its infinitesimal generator is denoted
by $-\calL$ and we call $\calL$ the Ornstein-Uhlenbeck operator
(sometimes called the number operator by the physicists). Due to the
Meyer inequalities (cf., for instance \cite{ASU}), the
norms defined by 
\begin{equation}
\label{norm}
\|\varphi\|_{p,k}=\|(I+\calL)^{k/2}\varphi\|_{L^p(\mu)}
\end{equation}
are equivalent to the norms defined by the iterates of the  Sobolev
derivative $\nabla$. This observation permits us to identify the duals
of the space $\DD_{p,k}(\Phi);p>1,\,k\in\NN$ by $\DD_{q,-k}(\Phi')$,
with $q^{-1}=1-p^{-1}$, 
where the latter  space is defined by replacing $k$ in (\ref{norm}) by
$-k$, this gives us the distribution spaces on the Wiener space $W$
(in fact we can take as $k$ any real number). An easy calculation 
shows that, formally, $\delta\circ \nabla=\calL$, and this permits us
to extend the  divergence and the derivative  operators to the
distributions as linear,  continuous operators. In fact
$\delta:\DD_{q,k}(H\otimes \Phi)\to \DD_{q,k-1}(\Phi)$ and 
$\nabla:\DD_{q,k}(\Phi)\to\DD_{q,k-1}(H\otimes \Phi)$ continuously, for
any $q>1$ and $k\in \reals$, where $H\otimes \Phi$ denotes the
completed Hilbert-Schmidt tensor product (cf., for instance \cite{ASU}). 
The following assertion is useful: assume that $(Z_n,n\geq
1)\subset \DD'$ converges to $Z$ in $\DD'$, assume further that each
each $Z_n$ is a probability measure on $W$, then $Z$ is also a
probability and $(Z_n,n\geq 1)$ converges to $Z$ in the weak topology
of measures. In particular, a lower bounded distribution (in the sense
that there exists a constant $c\in \R$ such that $Z+c$ is a positive
distribution) is a (Radon) measure on $W$, c.f. \cite{ASU}.

A  measurable  function
$f:W\to \reals\cup\{\infty\}$ is called $H$-convex (cf.\cite{F-U1}) if 
$$
h\to f(x+h)
$$
is convex $\mu$-almost surely, i.e., if for any $h,k\in H$, $s,t\in
[0,1],\,s+t=1$, we have 
$$
f(x+sh+tk)\leq sf(x+h)+tf(x+k)\,,
$$
almost surely, where the negligeable set on which this inequality
fails may depend on the choice of $s,h $ and of $k$. We can rephrase
this property by saying that $h\to(x\to f(x+h))$ is an $L^0(\mu)$-valued
convex function on $H$. 
$f$ is called $1$-convex if the map 
$$
h\to\left(x\to f(x+h)+\frac{1}{2}|h|_H^2\right)
$$
is convex on the Cameron-Martin space $H$ with values in
$L^0(\mu)$. Note that all these  notions are  compatible with the
$\mu$-equivalence classes of random variables thanks to the
Cameron-Martin theorem. It is proven in \cite{F-U1} that 
this definition  is equivalent  the following condition:
  Let $(\pi_n,n\geq 1)$ be a sequence of regular, finite dimensional,
  orthogonal projections of $H$,  increasing to the identity map
  $I_H$. Denote also  by $\pi_n$ its  continuous extension  to $W$ and
  define $\pi_n^\bot=I_W-\pi_n$. For $x\in W$, let $x_n=\pi_nx$ and
  $x_n^\bot=\pi_n^\bot x$.   Then $f$ is $1$-convex if and only if 
$$
x_n\to \frac{1}{2}|x_n|_H^2+f(x_n+x_n^\bot)
$$ 
is  $\pi_n^\bot\mu$-almost surely convex. 
We define similarly the notion of $H$-concave and $H$-log-concave
functions. In particular, one can prove that, for any $H$-log-concave
function $f$ on $W$, $P_tf$ and $E[f|V_n]$ are again $H$-log-concave
\cite{F-U1}. 

\section{Monge-Kantorovitch problem}

Let us recall the definition of the Monge-Kantorovitch problem in our
case:
\begin{definition}
Let $\rho$ and $\nu$ be two probability measures on $W$, let also
$\Sigma(\rho,\nu)$ be the convex subset of the probability measures on
the product space $W\times W$ whose first marginal is $\rho$ and the
second one is $\nu$. The Monge-Kantorovitch problem for the couple
$(\rho,\nu)$ consists of finding a measure $\ga\in \Sigma (\rho,\nu)$
which realizes the following infimum:
$$
d^2_H(\rho,\nu)=\inf_{\beta\in \Sigma(\rho,\nu)}\int_{W\times
  W}|x-y|_H^2d\beta(x,y)\,.
$$
The function $c(x,y)=|x-y|_H^2$ is called the cost function.
\end{definition}
\begin{remarkk}
Note that the cost function is not continuous with respect to the
product topology of $W\times W$ and it takes the value $\infty$ very
often for the most notable measures, e.g., when $\rho$ and $\nu$ are
absolutely continuous with respect to the Wiener measure $\mu$.
\end{remarkk}

\noindent
The proof of the next theorem, for which we refer the reader to
\cite{F-U3},  can be done by choosing a proper 
disintegration of any optimal measure in such a way that the elements
of this  disintegration are the solutions of finite dimensional
Monge-Kantorovitch problems. The latter is proven with the help of
the measurable section-selection theorem.

\begin{theorem}[General case]
\label{monge-general}
Suppose that $\rho$ and $\nu$ are two probability measures on
$W$ such  that
$$
d_H(\rho,\nu)<\infty\,.
$$
Let $(\pi_n,n\geq 1)$ be a total increasing  sequence of regular
projections (of $H$, converging to the identity map of $H$). 
Suppose  that, for any $n\geq 1$, the regular
conditional probabilities $\rho(\cdot\,|\pi_n^\bot=x^\bot)$ vanish
$\pi_n^\bot\rho$-almost surely on 
the subsets of  $(\pi_n^\bot)^{-1}(W)$ with Hausdorff dimension
$n-1$. Then there exists a    
unique solution of the   Monge-Kantorovitch problem, denoted by $\ga\in
\Sigma(\rho,\nu)$ and  $\ga$ is supported by the graph of a Borel
map $T$ which is the solution of the
Monge problem.  $T:W\to W$ is  of the form $T=I_W+\xi$ , where $\xi\in
H$ almost surely. Besides  we have 
\beaa
d_H^2(\rho,\nu)&=&\int_{W\times W}|T(x)-x|_H^2d\ga(x,y)\\
&=&\int_{W}|T(x)-x|_H^2d\rho(x)\,, 
\eeaa
and  
for $\pi_n^\bot\rho$-almost almost all $x_n^\bot$, the map $u\to
u+\xi(u+x_n^\bot)$ is cyclically monotone on
$(\pi_n^\bot)^{-1}\{x_n^\bot\}$, in the sense that 
$$
\sum_{i=1}^N\left(u_i+\xi(x_n^\bot+u_i),u_{i+1}-u_i\right)_H\leq 0
$$
$\pi_n^\bot\rho$-almost surely, for any cyclic sequence
$\{u_1,\ldots,u_N,u_{N+1}=u_1\}$ from $\pi_n(W)$. Finally, if, for any $n\geq
1$, $\pi_n^\bot\nu$-almost surely,  $\nu(\cdot\,|\pi_n^\bot=y^\bot)$
 also vanishes on the $n-1$-Hausdorff dimensional  subsets  
 of $(\pi_n^\bot)^{-1}(W)$, then $T$ is invertible, i.e, there exists
 $S:W\to W$ of the form $S=I_W+\eta$ such that  $\eta\in H$ satisfies
 a similar  cyclic monotononicity property as $\xi$ and that 
\beaa
1&=&\ga\left\{(x,y)\in W\times W: T\circ S(y)=y\right\}\\
&=&\ga\left\{(x,y)\in W\times W: S\circ T(x)=x\right\}\,.
\eeaa
In particular we have 
\beaa
d_H^2(\rho,\nu)&=&\int_{W\times W}|S(y)-y|_H^2d\ga(x,y)\\
&=&\int_{W}|S(y)-y|_H^2d\nu(y)\,. 
\eeaa
\end{theorem}
\begin{remarkk}
In particular, for all  the measures $\rho$ which are  absolutely
continuous with respect to the  Wiener measure $\mu$,  the second
hypothesis is satisfied, i.e., the measure
$\rho(\cdot\,|\pi_n^\bot=x_n^\bot)$ vanishes on the sets of Hausdorff
dimension $n-1$.
\end{remarkk} 

\noindent
Any probability measure satisfying the hypothesis of Theorem
\ref{monge-general} is called a spread measure. Namely, 
\begin{definition}
A probability measure $m$ on $(W,\calB(W))$ is called a spread measure if
there exists a sequence of finite dimensional regular  projections
$(\pi_n,n\geq 1)$ converging to $I_H$ such that the regular
conditional probabilities $m(\,\cdot\,|\pi_n^\bot=x_n^\bot)$ concentrated
in the $n$-dimensional spaces $\pi_n(W)+x_n^\bot$ vanishe on the sets
of Hausdorff dimension $n-1$ for $\pi_n^\bot(m)$-almost all
$x_n^\bot$ and for any $n\geq 1$. 
\end{definition}

The case where one of the measures is the Wiener measure  and the other is
absolutely continuous with respect to $\mu$ is the most important one
for the applications. Consequently we give the related results
separately in the following theorem where the tools of the Malliavin
calculus give more information about the maps $\xi$ and $\eta$ of
Theorem \ref{monge-general}:  
\begin{theorem}[Gaussian case]
\label{gaussian-case}
Let $\nu$ be the measure $d\nu=Ld\mu$, where $L$ is a positive random variable,
with $E[L]=1$. Assume that $d_H(\mu,\nu)<\infty$ (for instance 
$L\in \LL\log \LL$). Then there exists a  $1$-convex function $\varphi\in
\DD_{2,1}$ and  a partially $1$-convex function $\psi\in L^2(\nu)$,
both are unique upto a constant, called Monge-Kantorovitch potentials,
such that   
$$
\varphi(x)+\psi(y)+\frac{1}{2}|x-y|_H^2\geq 0
$$
for all $(x,y)\in W\times W$ and that 
$$
\varphi(x)+\psi(y)+\frac{1}{2}|x-y|_H^2= 0
$$
$\ga$-almost everywhere. 
The map $T=I_W+\nabla \varphi$ is the unique solution of the original problem of
Monge. Moreover, its graph supports  the  unique
solution of the 
Monge-Kantorovitch problem $\ga$. Consequently    
$$
(I_W\times T)\mu=\ga
$$
In particular  $T$ maps $\mu$ to $\nu$ and  $T$ is almost surely
invertible, i.e., there exists some $T^{-1}=I_W+\eta$ such that
$T^{-1}\nu=\mu$, $\eta\in L^2(\nu)$
and that 
\beaa
1&=&\mu\left\{x:\,T^{-1}\circ T(x)=x\right\}\\
&=&\nu\left\{y\in W:\,T\circ T^{-1}(y)=y\right\}\,.
\eeaa
\end{theorem}

\begin{remarkk}
\label{nu-closed}
By the partial $1$-convexity we mean that $y_F\to \psi(y_F+y^\bot_F)$
is $\nu(\cdot|\pi_F^\bot=y_F^\bot)$-almost surely  $1$-convex on any regular, finite dimensional
subspace $F$, where $\pi_F$ is the (regular) projection corresponding
to $F$, $y_F=\pi_F(y)$ and  $y_F^\bot=y-y_F$. Assume that the operator
$\nabla$ is closable with respect to $\nu$, then we have
$\eta=\nabla\psi$. In particular, if $\nu$ and $\mu$ are equivalent,
then we have    
$$
T^{-1}=I_W+\nabla\psi\,,
$$
where is $\psi$ is  a $1$-convex function.
\end{remarkk}
\label{stability}
\begin{remarkk}
Let $(e_n,n\in \NN)$ be a complete,  orthonormal in $H$, denote by
$V_n$ the sigma algebra generated by $\{\delta e_1,\ldots,\delta
e_n\}$ and let  $L_n=E[L|V_n]$. If $\varphi_n\in \DD_{2,1}$ is the
function constructed in Theorem \ref{gaussian-case}, corresponding to
$L_n$, then, using the inequality (cf., \cite{F-U3})
$$
d_H^2(\mu,\nu)\leq 2E[L\log L]\,,
$$
we can prove that
the sequence $(\varphi_n,n\in \NN)$ converges to $\varphi$ in $\DD_{2,1}$. 
\end{remarkk}

\section{Polar factorization of mappings between  spread measures}
\label{pol-fac}
In \cite{F-U3} we have proved the polar factorization of the mappings
$U:W\to W$ such that the Wasserstein distance between $U\mu$ and the
Wiener measure $\mu$, denoted by  $d_H(\mu,U\mu)$,  is finite. We have
also studied the particular case where $U$ is a 
perturbations of identity, i.e., 
it is  the form $I_W+u$, where $u$ maps $W$ to the Cameron-Martin space
$H$.  In this section we shall
generalize this results  in the frame of spread measures.
\begin{theorem}
\label{polar-fac-thm}
Assume that $\rho$ and $\nu$ are spread measures with
$d_H(\rho,\nu)<\infty$  and that $U\rho=\nu$,
for some measurable map $U:W\to W$. Let $T$ be the optimal transport map
sending $\rho$ to $\nu$, whose existence and uniqueness is proven in
Theorem \ref{monge-general}.  Then $R=T^{-1}\circ
U$ is a $\rho$-rotation (i.e., $R\rho=\rho$) and $U=T\circ R$,
morover, if $U$ is a perturbation of identity, then $R$ is also a
perturbation of identity. In both cases, $R$  is the $\rho$-almost everywhere
unique  minimal $\rho$-rotation in the sense that  
\begin{equation}
\label{min-eq}
\int_W|U(x)-R(x)|_H^2d\rho(x)=\inf_{R'\in
  \mathcal{R}}\int_W|U(x)-R'(x)|_H^2d\rho(x)\,,
\end{equation}
where $\mathcal{R}$ denotes the set of $\rho$-rotations.
\end{theorem}
\proof
Let $T$ be the optimal transportation of $\rho$ to $\nu$ whose
existence and uniqueness follows from  Theorem \ref{monge-general}. The
unique solution $\ga$ of the Monge-Kantorovitch problem for
$\Sigma(\rho,\nu)$ can be written as $\ga=(I\times T)\rho$. Since
$\nu$ is spread, $T$ is invertible on the support of $\nu$ and we have
also $\ga=(T^{-1}\times I)\nu$. In particular $R\rho=T^{-1}\circ
U\rho=T^{-1}\nu=\rho$, hence $R$ is a rotation.   Let $R'$
be another rotation in $\mathcal R$, define $\ga'=( R'\times U)\rho$,
then $\ga'\in \Sigma(\rho,\nu)$ and the optimality of $\ga$ implies
that $J(\ga)\leq J(\ga')$, besides  we have 
\beaa
\int_W|U(x)-R(x)|_H^2d\rho(x)&=&\int_W|U(x)-T^{-1}\circ
U(x)|_H^2d\rho(x)\\
&=&\int_W|x-T^{-1}(x)|_H^2d\nu(x)\\
&=&\int_W|T(x)-x|_H^2d\rho(x)\\
&=&J(\ga)\,.
\eeaa
On the other hand
$$
J(\ga')=\int_W|U(x)-R'(x)|_H^2d\rho(x)\,,
$$
hence the relation (\ref{min-eq}) follows. Assume now that for  the
second rotation $R'\in \mathcal{R}$ we have the equality 
$$
\int_W|U(x)-R(x)|_H^2d\rho(x)=\int_W|U(x)-R'(x)|_H^2d\rho(x)\,.
$$
Then we have $J(\ga)=J(\ga')$, where $\ga'$ is defined above. By the
uniqueness of the solution of Monge-Kantorovitch problem due to Theorem
\ref{monge-general}, we should 
have $\ga=\ga'$. Hence $(R\times U)\rho=(R'\times U)\rho=\ga$,
consequently, we have 
$$
\int_W f(R(x),U(x))d\rho(x)=\int_W f(R'(x),U(x))d\rho(x)\,,
$$
for any bounded, measurable map $f$ on $W\times W$. This implies in
particular
$$
\int_W (a\circ T\circ R)\,\,(b\circ U) d\rho=\int_W (a\circ T\circ
R')\,\,(b\circ U) d\rho 
$$ 
for any bounded measurable functions $a$ and $b$.
Let $U'=T\circ R'$, then the above expression reads as 
$$
\int_W a\circ U\,\,b\circ U d\rho=\int_W a\circ U'\,\,b\circ U
d\rho\,.
$$
Taking $a=b$, we obtain 
$$
\int_W(a\circ U) \,(a\circ U')\,d\rho=\|a\circ U\|_{L^2(\rho)}\|a\circ
U'\|_{L^2(\rho)}\,,
$$
for any bounded, measurable $a$. This implies that $a\circ U=a\circ
U'$ $\rho$-almost surely for any $a$, hence $U=U'$  i.e, $T\circ
R=T\circ R'$$\rho$-almost surely. Let us denote by $S$ the left
inverse of $T$ whose existence follows from Theorem
\ref{monge-general} and let  $D=\{x\in W: S\circ T(x)=x\}$. Since
$\rho(D)=1$ and since $R$ and $R'$ are $\rho$-rotations,  we have also 
$$
\rho\left(D\cap R^{-1}(D)\cap R'^{-1}(D)\right)=1\,.
$$
Let  $x\in W$ be  any element of $ D\cap R^{-1}(D)\cap R'^{-1}(D)$, then 
\beaa
R(x)&=&S\circ T\circ R(x)\\
&=&S\circ T\circ R'(x)\\
&=&R'(x)\,,
\eeaa
consequently $R=R'$ on a set of full  $\rho$-measure..
\qed

Let us give another result of interest as an application of these
factorization results: it is important to have as much as information
about the measures and the tranformations which induce them in the
setting of Girsanov Theorem, cf. \cite{BOOK} and the references
there. The problem which we propose is the following: assume that,  in
the case of the Wiener measure,  we have a density $L$ with
$d_H(\mu,L\cdot\mu)<\infty$, hence from Theorem \ref{monge-general}, a
map $T:W\to W$, which is the optimal transport map corresponding to
the solution of MKP in $\Sigma(\mu,L\cdot\mu)$. Since the target
measure is also spread, the map $T$ possesses a left inverse $S$ such
that $S\circ T=I_W$ $\mu$-almost surely. Assume now that the
transformation $T$ has a Girsanov density, i.e., $\la\in L^1_+(\mu)$,
with $E[\la]=1$ and that 
$$
\int f\circ T \la \,d\mu=\int f\,d\mu\,,
$$
for any $f\in C_b(W)$. We can now prove:
\begin{theorem}
\label{girsanov-complement}
Let $T$ be as explained above, assume moreover that
$$
d_H(\la\cdot\mu,\mu)<\infty\,,
$$ 
then $T$ has also a right inverse,
i.e., $T$ is invertible $\mu$-almost everywhere.
\end{theorem}
\proof
Denote by $\Theta:W\to W$ the optimal transportation map corresponding
to the solution of MKP in $\Sigma(\mu,\la\cdot\mu)$. Note that both of
the measures  $(T\times I_W)(\la\cdot\mu)$ and $(I_W\times \Theta)\mu$
belong to $\Sigma(\mu,\la\cdot\mu)$.  By the uniqueness
of the solutions of MKP, they are equal, hence, for any $a,b\in
C_b(W)$, we have 
\begin{equation}
\label{iden-1}
\int a(T(x))\,b(x)\la(x)d\mu(x)=\int a(x)b(\Theta(x))d\mu(x)\,.
\end{equation}
Since $\Theta(\mu)=\la\cdot\mu$,  the equality
(\ref{iden-1}) can  also be written as 
\begin{equation}
\label{iden-2}
\int a(T\circ\Theta(x))\,b(\Theta(x))d\mu(x)=\int a(x)b(\Theta(x))d\mu(x)\,.
\end{equation}
Since, as $T$, the map $\Theta$ has also a left inverse, the sigma
algebra generated by $\Theta$ is equal to the Borel sigma algebra of
$W$, consequently, the relation (\ref{iden-2}) implies that 
$$
a\circ T\circ\Theta=a\,,
$$
$\mu$-almost surely, for any $a\in C_b(W)$. Therefore we have 
$$
\mu\left(\{x\in W:\,T\circ \Theta(x)=x\}\right)=1\,,
$$
since $T$ has already a left inverse, the proof is completed.
\qed

\subsection{Application to Gaussian measures}
\label{example-1}
{\rm{Let us give an example of the above results: Assume that $\rho=\mu$,
i.e., the Wiener measure and let $K$ be a Hilbert-Schmidt operator on
$H$. Assume that the Carleman-Fredholm determinant $\dett(I_H+K)$ is
different than zero, hence  the operator  $I_H+K:H\to H$ is
invertible. Moreover, it follows from the general theory that $I_H+K$
has a unique polar decomposition as $I_H+K=(I_H+\bar{K})(I_H+A)$,
where $I_H+A$ is an isometry{\footnote{$A$ satisfies the relation
    $A+A^*+A^*A=0$.}} and $I_H+\bar{K}$ is a symmetric, 
positive operator. Note that $\bar{K}$ is compulsorily
Hilbert-Schmidt. Let us now define $U:W\to W$ as $U(x)=x+\delta K(x)$,
where $\delta K(x)$ is the $H$-valued divergence, defined by $(\delta
K(x),h)_H=\delta(K^*h)(x)$. Then it is known that the measure $U\mu$
is absolutely continuous with respect to $\mu$, in fact $U\mu$ is even
equivalent to $\mu$ since $|\La_K|\neq 0$ $\mu$-almost surely, where 
$$
\La_K=\dett(I_H+K)\exp\left\{\delta^2(K)-\frac{1}{2}|\delta
  K|_H^2\right\}\,.
$$
Besides we have 
$$
L=\frac{dU\mu}{d\mu}=\frac{1}{|\La_K|\circ V}\,,
$$
where $V$ is the inverse of $U$, whose existence follows from the
invertibility of $h\to h+\delta(K)(x)+Kh$ on H, cf.
\cite{BOOK}. Consequently, 
$$
E[L\log L]=-E[\log |\La_K|]<\infty\,,
$$
hence $d_H(\mu,U\mu)<\infty$. We shall prove that the polar
factorization of $U$ is given by 
$$
U= (I_W+\delta\bar{K})\circ(I_W+\delta A)\,.
$$
In fact, it follows from Theorem B.6.4 of \cite{BOOK}, that 
\beaa
 (I_W+\delta\bar{K})\circ (I_W+\delta A)&=&I_W+\delta\bar{K}+\delta
A+\delta(\bar{K}A)\\ 
&=&I_W+\delta(\bar{K}+A+\bar{K}A)\\
&=&I_W+\delta K\,.
\eeaa
Besides $\nabla^2\delta^2\bar{K}=2\bar{K}$, and since $I_H+\bar{K}$ is
a positive operator, the Wiener map $\frac{1}{2}\delta^2\bar{K}$ is $1$-convex,
consequently, $T=I_W+\delta\bar{K}$ is the transport map and
$I_W+\delta A$ is the unique rotation whose existence is proven in
Theorem \ref{polar-fac-thm}. The Kantorovitch potentials $\varphi$ and
$\psi$ of Theorem \ref{gaussian-case} can be chosen as 
$$
\varphi(x)=\frac{1}{2}\delta^2\bar{K}(x)
$$ 
for $T$ and
$$
\psi(x)=-\frac{1}{2}\delta((I_H+\bar{K})^{-1}\bar{K})(x)
$$ 
for
$T^{-1}=I_W+\nabla \psi$.
}}
\begin{remarkk}
{\rm{ Let us denote by $P_{\ker\delta}$ the projection operator from
    $\DD'(H)$ to the kernel of the divergence operator $\delta$. Then,
    we have the following identity:
$$
P_{\ker\delta}\left(\delta((I_H+\bar{K})A)\right)=
\delta\hat{K}-\delta\bar{K}\,,
$$
where $\hat{K}$ denotes the symmetrization of $K$. This shows that the
polar decomposition and the Helmholtz decomposition are different in
general. 

}}
\end{remarkk}
We can also calculate the Monge-Kantorovitch potential function for
the singular case as follows: assume that $\nu$ is a zero mean
Gaussian measure on $W$ such that $d_H(\mu,\nu)<\infty$. Then there
exists a bilinear form $q$ on $W^\star$ such that 
$$
\int_W
e^{i\langle\alpha,x\rangle}d\nu(x)=\exp-\frac{1}{2}q(\alpha,\alpha)\,,
$$
for any $\alpha\in W^\star$. On the other hand, from Theorem
\ref{gaussian-case}, there exists a $\varphi\in\DD_{2,1}$, which is
$1$-convex, such that $T\mu=(I_W+\nabla\varphi)\mu=\nu$. Hence, rewriting
the above relation with $T$, we obtain:
\begin{equation}
\label{car-func}
\int_W
e^{i\langle t\alpha,T(x)\rangle}d\mu(x)=\exp-\frac{t^2}{2}q(\alpha,\alpha)\,,
\end{equation}
for any $t\in \R$ and $\alpha\in W^\star$. Taking the derivative of
both sides twice  at $t=0$, we obtain 
\begin{eqnarray*}
q(\alpha,\alpha)&=&|\tilde{\alpha}|_H^2+E\left[(\nabla\varphi,\tilde{\alpha})_H^2\right]+
2E\left[(\nabla\varphi,\tilde{\alpha})_H\delta\tilde{\alpha}\right]\\
&=&|\tilde{\alpha}|_H^2+E\left[(\nabla\varphi\otimes\nabla\varphi,\tilde{\alpha}\otimes\tilde{\alpha})_2\right]+2E\left[(\nabla^2\varphi,\tilde{\alpha}\otimes\tilde{\alpha})_2\right]\,, 
\end{eqnarray*}
where $\tilde{\alpha}$ denotes the image of $\alpha$ under the
injection $W^\star\hookrightarrow H$.Note that, here,  $\nabla^2\varphi$ is to be
interpreted as a distribution. Denote by $M$ the Hilbert-Schmidt
operator defined by 
$$
M=E\left[\nabla\varphi\otimes\nabla\varphi\right]+2E\left[\nabla^2\varphi\right]\,.  
$$
We have 
$$
q(\alpha,\alpha)=((I_H+M)\tilde{\alpha},\tilde{\alpha})_H\,.
$$
Let $I_H+N$ be the positive  square root of the (positive) operator
$I_H+M$, then $N$ is a symmetric Hilbert-Schmidt operator. Define 
$$
\varphi=\frac{1}{2}\delta^2N\,.
$$ 
Evidently $\varphi$ is a $1$-convex element of $\DD_{2,1}$, moreover
the map $T$ defined by $T=I_W+\nabla\varphi=I_W+\delta N$ satisfies
the identity (\ref{car-func}), hence $T$ is the unique solution of
the Monge problem and $(I_W\times T)\mu$ is the unique solution of MKP
for $\Sigma(\mu,\nu)$.

\section{Strong solutions  of the  Monge-Amp\`ere equation for
  $H$-log-concave  densities} 
\label{st-MA}
Assume that $L\in \LL^1_{+,1}(\mu)$ is of the form
$$
L=\frac{1}{E\left[e^{-f}\right]}e^{-f}\,,
$$
where $f$ is an  $H$-convex function in some $L^p(\mu)$, $p>1$. We
assume that $f\geq -\alpha$ almost surely, for some  $\alpha\in \R_+$. 
Denote by $\varphi\in \DD_{2,1}$ the potential of
the transport problem between $\mu$ and $\nu=L\cdot \mu$ which is a
$1$-convex function. This means that
the mapping defined by  $T=I_W+\nabla\varphi$ satisfies 
$T\mu=L\cdot\mu$ and $(I_W\times T)\mu$ is the unique  solution of the
Monge-Kantorovitch problem in $\Sigma(\mu,\nu)$ with the
singular quadratic cost function 
$c(x,y)=|x-y|_H^2$. Let $\La=1/L\circ T$, we know that 
$T^{-1}\mu=\La\cdot \mu$ where $T^{-1}=I_W+\nabla\psi$ such that
$\psi\in L^2(\nu),\,\nabla\psi\in L^2(\nu,H)$ (cf. Remark
\ref{nu-closed}) is also defined  uniquely\footnote{In fact in the
  proof of Theorem \ref{D-22-thm}, we shall see that
  $(\psi_n,n\geq 1)$ is bounded in $\DD_{2,1}$.}. Let $L_n=E[L|V_n]$, where  
$V_n$ is the sigma 
algebra generated by the first $n$ elements of an orthonormal basis
$(e_n,n\geq 1)$ of $H$. It follows from \cite{F-U1}, that $L_n$ is of
the form $\frac{1}{c} e^{-f_n}$, where $f_n$ is an $H$-convex function
on $W$ and $c=E[e^{-f}]$. 
We denote by $\varphi_n,\,\La_n,\,\psi_n$ the maps associated to
$L_n$, i.e., $T_n=I_W+\nabla\varphi_n$ maps $\mu$ to the measure
$L_n\cdot \mu$ and $S_n=I_W+\nabla\psi_n$ maps $L_n\cdot\mu$ to $\mu$.
Besides, from \cite{Caf}, $\nabla\varphi_n$ is a $1$-Lipschitz
map, i.e., 
$$
|\nabla\varphi_n(x)-\nabla\varphi_n(y)|\leq |x-y|\,,
$$
for any $x,y\in \R^n$, here it is  remarkable that  the Lipschitz
constant  is one and it is independent of the dimension of the
underlying space. Hence $\calL \varphi_n$ is a well-defined
element of $L^2(\mu)$, 
$|\nabla\varphi|_H^2$ is exponentially integrable, i.e., there exists
some $t>0$ such that 
\begin{equation}
\label{exp-bound}
\sup_nE\left[\exp t|\nabla\varphi_n|_H^2\right]<\infty\,,
\end{equation}
then the Fatou Lemma implies that 
$$
E\left[\exp t|\nabla\varphi|_H^2\right]<\infty\,.
$$
It follows in particular  that $(\varphi_n,n\geq 1)\subset  \DD_{p,2}$
and it  converges to $\varphi$ in $\DD_{p,1}$ for any
$p\geq 1$, cf., \cite{F-U3}. Moreover, from a result of McCann
\cite{Mc2}, we  have  
$$
\La_n=\dett(I_H+\nabla^2\varphi_n)\exp\left\{-\calL\varphi_n-\frac{1}{2}
|\nabla\varphi_n|_H^2\right\}\,.
$$
Since $\La_n=1/L_n\circ T_n$, the sequence $(\La_n,n\geq 1)$ is lower 
bounded. Hence $(-\log \La_n, n\geq 1)$ is upper  bounded, besides 
\beaa
E[|\log \La_n|^p]&\leq &C_p E[|f_n\circ T_n|^p]+D_p\\
&=&C_pE[|f_n|^pL_n]+D_p\\
&\leq&C_p e^{p\alpha}E[|f|^p]+D_p\,,
\eeaa
where $C_p$ and $D_p$ are some constants. Since $(-\log\La_n,n\geq 1)$
converges in $\LL^0(\mu)$ to $-\log\La$, it follows from the dominated
convergence theorem that $(\log \La_n,n\geq 1)$ converges to $\log\La$ in
$\LL^p(\mu)$. Therefore 
$$
-\log\dett(I_H+\nabla^2\varphi_n)+\calL\varphi_n+\frac{1}{2}
|\nabla\varphi_n|_H^2\to-\log\La
$$
in $\LL^p(\mu)$. Since $(\varphi_n,n\geq 1)$ converges
to $\varphi$ in $\cap_p \DD_{p,1}$, the sequence $(Z_n,n\geq 1)$, defined by 
$$
Z_n=-\log\dett(I_H+\nabla^2\varphi_n)+\calL\varphi_n\,,
$$
converges in $\LL^p(\mu)$ to some $Z\in \LL^p(\mu)$. Again by the
convergence of $(\varphi_n,n\geq 1)$, the sequence $(\calL\varphi_n,n\geq
1)$ converges to the measure $\calL\varphi$ in $\DD_{2,-1}$
(cf. \cite{F-U3}), consequently 
the sequence $(\log\dett(I_H+\nabla^2\varphi_n),n\geq 1)$ converges to
some $D=D(\varphi)$ in $\DD'$. Since $Z=\calL\varphi-D(\varphi)$ and $\calL\varphi$ are  measures,
$D(\varphi)$ should be a measure, besides $Z$ is 
absolutely continuous with respect to $\mu$ (it is a random variable),
hence $\calL_s\varphi-D_s(\varphi)=0$, where the subscript ``s'' denotes the
singular part of the  measure $D(\varphi)$. Consequently we have
$Z=\calL_a\varphi-D_a(\varphi)$, where the subscript ``$a$'' denotes
the absolutely continuous part of the corresponding measure. Therefore
we have  
\beaa
\La&=&\lim\La_n\\
&=&\exp\left\{D_a(\varphi)-\calL_a\varphi-\half|\nabla\varphi|_H^2\right\}\,.
\eeaa

\noindent
In fact we have a much better result of regularity:
\begin{theorem}
\label{D-22-thm}
Assume further  that $f\in \DD_{2,1}$, then $\varphi\in \DD_{2,2}$, in
particular 
$$
\calL_a\varphi=\calL\varphi\in L^2(\mu)
$$ 
and
$\dett(I_H+\nabla^2\varphi)$ is a  well-defined function.
\end{theorem}

In order to  proceed to the proof of Theorem \ref{D-22-thm}, we need a
lemma whose proof is given in a more general case in \cite{BOOK},
Appendix B:
\begin{lemma}
\label{commutation-lemma}
Assume that $M:W\to W$ is a map of the form $M=I_W+u$, where
$u\in\DD_{2,1}(H)$ such that $M\mu$ is absolutely continuous with
respect to $\mu$. For any smooth, cylindrical vector field $\xi:W\to
H$, we have 
$$
\delta\xi\circ M=\delta(\xi\circ M)+(\xi\circ
M,u)_H+\trace(\nabla\xi\circ M\cdot\nabla u)\,,
$$
$\mu$-almost surely.
\end{lemma}
\proof It suffices to represent $\xi$ with an orthonormal basis
$(e_i,i\geq 1)$ of $H$ as 
$$
\xi=\sum_i(\xi,e_i)e_i\,,
$$
then
$$
\delta(\xi\circ M)=\sum_i(\xi\circ M,e_i)_H\delta
e_i-\nabla_{e_i}(\xi\circ M,e_i)_H\,.
$$
Since $\delta e_i\circ M=\delta e_i+(e_i,\xi)_H$ and since
$\nabla(\xi\circ M)=\nabla\xi\circ M(I_H+\nabla\xi)$, we obtain at
once the claimed equality.
\qed

\noindent
{\bf{Proof of Theorem \ref{D-22-thm}:}}
 $L_n$ is $\mu$-a.s. strictly positive by the hypothesis that we have
 done for  $L$. Consequently, the operator $I+\nabla^2\varphi_n(x)$ is
 almost surely invertible. Besides, using Lemma \ref{commutation-lemma}
 and the relation $\delta\circ \nabla=\calL$, we get 
\begin{equation}
\label{key-com}
\calL\psi_n\circ T_n=\delta(\nabla\psi_n\circ T_n)+(\nabla\psi_n\circ
T_n,\nabla\varphi_n)_H+\trace(\nabla^2\psi_n\circ
T_n\,\cdot\nabla^2\varphi_n)\,.
\end{equation}
It is easy to see that 
$$
\trace(\nabla^2\psi_n\circ
T_n\,\cdot\nabla^2\varphi_n)=-\trace\left((I+\nabla^2\varphi_n)^{-1}\cdot(\nabla^2\varphi_n)^2\right)\,.
$$
Taking the expectation of both sides of (\ref{key-com}) with respect
to $\mu$, we have 
$$
E\left[\trace\left((I+\nabla^2\varphi_n)^{-1}\cdot(\nabla^2\varphi_n)^2\right)\right]
=E[|\nabla\varphi_n|_H^2]-E[\calL\psi_n\,L_n]\,.
$$
Since $(L_n,n\geq 1)$ is uniformly essentially bounded by some $K>0$, we have 
\beaa
E[\calL\psi_n\,L_n]&=&E[(\nabla\psi_n,\nabla L_n)_H]\\
&=&-E[(\nabla\psi_n,\nabla f_n)_H L_n]\\
&\leq& K\,\|\nabla\psi_n\|_{L^2(\mu,H)}\|f\|_{2,1}\,.
\eeaa
Moreover, from the Young inequality
$$
E[|\nabla\psi_n|_H^2]=E[|\nabla\varphi_n|_H^2\,\La_n]\leq
E[\eps^{-1}\La_n\log\La_n]+E\left[\exp\eps|\nabla\varphi_n|_H^2\right]\,, 
$$
which is uniformly bounded with respect to $n$ since
$\|\nabla^2\varphi_n\|_{\rm{op}}\leq 1$. 
Consequently
$$ 
\sup_nE\left[\trace\left((I+\nabla^2\varphi_n)^{-1}\cdot(\nabla^2\varphi_n)^2\right)\right]<\infty\,.
$$
Recalling that $\|I_H+\nabla^2\varphi_n\|_{\rm{op}}\leq 1$ almost
surely, we finally get 
\beaa
\sup_n E\left[\trace(\nabla^2\varphi_n)^2\right]&=&
\sup_nE[\|\nabla^2\varphi_n\|_2^2]\\
&\leq&\sup_nE\left[\left\|(I_H+\nabla^2\varphi_n)^{-1/2}\nabla^2\varphi_n\right\|_2^2\right]\\
&=&
\sup_nE\left[\trace\left((I+\nabla^2\varphi_n)^{-1}\cdot(\nabla^2\varphi_n)^2\right)\right]<\infty\,.
\eeaa
This implies that $(\nabla^2\varphi_n,n\geq 1)$ is bounded in the
space Hilbert-Schmidt valued Wiener maps  $L^2(\mu,H\otimes H)$, since
$(\varphi_n, n\geq 1)$ converges to $\varphi$ in $\DD_{2,1}$, $\varphi$
should be in $\DD_{2,2}$ and the other claims are now immediate.
\qed

\begin{corollary}
\label{subsoln}
Let $\la$ be the function defined as 
$$
\la=\dett(I_H+\nabla^2\varphi)\exp\left\{-\calL\varphi-\frac{1}{2}|\nabla\varphi|_H^2\right\}\,.
$$
Then $\la$ is a sub-solution of the Monge-Amp\`ere equation in the
sense that 
\begin{equation}
\label{ineq-10}
E[g\circ T\,\la]\leq E[g]\,,
\end{equation}
for any positive, measurable function $g$. In particular 
$$
\la\leq \La
$$
almost surely.
\end{corollary} 
\proof
Let $(e_n,n\geq 1)\subset W^*$ be a complete, orthonormal basis of
$H$, denote by $V_n$ the sigma algebra generated by $\{\delta
e_1,\ldots,\delta e_n\}$. Since, from Theorem \ref{D-22-thm},
$\varphi\in \DD_{2,2}$, the sequence 
$(F_n,n\geq 1)$, where $F_n=E[\varphi|V_n]$ , converges to $\varphi$ in
$\DD_{2,2}$, hence the sequence $(M_n,n\geq 1)$, where
$$
M_n=\dett(I_H+\nabla^2F_n)\exp\left\{-\calL F_n-\frac{1}{2}|\nabla
F_n|_H^2\right\}\,, 
$$
converges to $\la$ in probability. Since $F_n$ is a $1$-convex
function, it follows from Theorem 6.3.1 of \cite{BOOK} that
$$
E[g\circ(I_W+\nabla F_n)\,M_n]\leq E[g]\,,
$$
for any positive, measurable function $g$.
The first  claim  follows  from the Fatou
lemma. Since $L>0$ almost surely, we have 
\begin{equation}
\label{ineq-11}
E[g\circ T\,\La]=E[g]\,,
\end{equation}
for any positive, measurable $g$, where 
$$
\La=\frac{1}{L\circ T}\,.
$$
As $T$ is invertible, we get $\la\leq \La$ by comparing the
relations (\ref{ineq-10}) and (\ref{ineq-11}).

\qed

\noindent
We can prove now the main theorem of this section:

\begin{theorem}
\label{jacobian-calcul}
Let $L$ be given as $c^{-1}\,e^{-f}$, where $f\in\DD_{2,1}$ is a lower
bounded, finite, 
$H$-convex Wiener function and define the probability  measure $\nu$
as $d\nu=Ld\mu$, where $c=E[e^{-f}]$ is the  normalization constant. Let
$T=I_W+\nabla \varphi$ be the
optimal transportation of $\mu$ to $\nu$ in the sense of Wasserstein
distance, where $\varphi\in \DD_{2,1}$ is the  $1$-convex potential
function. Then $\varphi\in \DD_{2,2}$ and  the Gaussian
Jacobian  of $T$ is equal to $\La=1/L\circ T$ and we have the
following relation: 
\be
\label{jacobian}
\La=\dett(I_H+\nabla^2\varphi)\exp\left\{-\calL\varphi-
\half|\nabla\varphi|_H^2\right\}\,.    
\ee
\end{theorem}
\proof
We have prepared everything necessary for the proof. First, we can
form a sequence, denoted by $\varphi_n',n\geq 1)$ such that each
$\varphi'_n$ is obtained as a convex combination from the elements
of the  tail sequence $(\varphi_k,k\geq n)$ and  that the sequence
$(\varphi'_n,n\geq 1)$  converges to $\varphi$ in $\DD_{2,2}$. Let us
denote the Jacobian  written with $\varphi_n'$  by $\La_n(\varphi_n')$
whose explicit expression is  given as 
$$
\La_n(\varphi_n')=\dett(I+\nabla^2\varphi'_n)\exp\left\{-\calL
\varphi'_n-\frac{1}{2}|\nabla\varphi'_n|_H^2\right\}
$$
Let $T'_n=I_W+\nabla\varphi'_n$ and $S'_n=I_W+\nabla\psi'_n$. Since
$A\to -\log\dett(I_H+A)$ is a convex function on the space of
symmetric Hilbert-Schmidt operators which are lower bounded by $-I_H$
(cf. \cite{B-B}, p.63), we have 
\beaa
-\log\La_n(\varphi_n')&=&-\log\dett\left(I_H+\sum_it_i\nabla^2\varphi_{n_i}\right)\\
&&+\sum_it_i\calL\varphi_{n_i}+\frac{1}{2}\left|\sum_it_i\nabla\varphi_{n_i}\right|_H^2\\
&\leq&\sum_i-t_i\log\La_{n_i}\,.
\eeaa
Since $(-\log \La_n,n\geq 1)$ converges to $-\log\La$ in any $L^p$
and since $(-\log\La_n(\varphi_n'),n\geq 1)$ converges to 
$-\log\la$,  it follows from the above inequality that 
$$
-\log\la\leq -\log\La
$$
almost surely, consequently $\La\leq \la$ almost surely. It follows
then from Corollary \ref{subsoln} that $\la=\La$ almost surely and
this completes the proof.
\qed

The following corollary gives the exact value of the Wasserstein
distance:
\begin{corollary}
\label{distance-cor}
With the hypothesis of Theorem \ref{jacobian-calcul}, we have 
$$
\frac{1}{2}d_H^2(\mu,L\cdot\mu)=E[L\log
L]+E\left[\log\dett(I_H+\nabla^2\varphi)\right]\,.
$$
\end{corollary}
\proof
Since $\La=c\,e^{f\circ T}$, it follows from the theorem that 
\beaa
\frac{1}{2}d_H^2(\mu,L\cdot\mu)&=&-E[f\circ T]-\log
c+E\left[\log\dett(I_H+\nabla^2\varphi)\right] \\
&=&E[L\log L]+E\left[\log\dett(I_H+\nabla^2\varphi)\right]\,.
\eeaa
In particular, the fact that
$E\left[\log\dett(I_H+\nabla^2\varphi)\right]$ is always negative
explains the defect in the  Talagrand inequality \cite{Tal}.
\qed

\noindent

Let us give an interesting result about the upper bound of the
interpolated density whose proof makes use also the convexity results
as in the proof of Theorem \ref{jacobian-calcul} :

\begin{proposition}
\label{inter-bound}
Assume the hypothesis of Theorem \ref{jacobian-calcul}, in particular
the relation $f\geq -\alpha$.  Denote by
$T_t=I_W+t\nabla \varphi$, $t\in [0,1]$,  then the Radon-Nikodym
density  $L_t$ the measure $T_t\mu $ with respect to $\mu$, 
satisfies the following inequality:
$$
L_t\leq \frac{1}{c}\exp\alpha t
$$
almost surely, where $c=E[\exp-f]$.
\end{proposition}
\proof
Let $g$ be any positive, measurable function on $W$, by the
convexity of $t\to -\log\La_t$, we have $-\log \La_t\leq -t\log\La$. Therefore
\beaa
E[L_t\log L_t\,g]&=&E[-\log \La_t\,g\circ T_t]\\
&\leq&E[-t\log\La\,g\circ T_t]\\
&=&E[-t(f\circ T+\log c)g\circ T_t]\\
&\leq&E[(t\alpha-\log c)g\circ T_t]\\
&=&E[(t\alpha-\log c)L_t\,g]\,.
\eeaa
Consequently 
$$
L_t\log L_t\leq (t\alpha-\log c)\,L_t
$$
almost surely.
\qed


\section{ It\^o-solutions of the Monge-Amp\`ere equation}
\label{Ito-section}
In the following calculations we shall take  $W$ as  the classical
Wiener space $W=C_0([0,1],\R)$, $H=H^1$, i.e., the Sobolev space
$W_{2,1}([0,1])$. We note that this choice does not entail any
restriction of generality as indicated in \cite{BOOK}, Chapter 2.6. 
Suppose we are given a positive  random variables $L=\frac{1}{c}e^{-f}$  whose
expectation is  equal to one, $c$ being the normalization
constant. Define the measure $\nu$ as  $d\nu=Ld\mu$. We shall suppose
that the Wasserstein 
distance $d_H(\mu,\nu)$ is finite, hence the conclusions of Theorem
\ref{monge-general} are valid.  In order to simplify the discussion we
shall assume that  $L$ is  strictly positive. The transport map $T$
can be represented as 
$T=I_W+\nabla\varphi$ again with $\varphi\in \DD_{2,1}$.  Define now 
$$
\La=\frac{1}{L\circ T}\,.
$$
We have 
$$
\int g\circ T\,\,\La\,d\mu=\int g\,d\mu\,,
$$
for any $g\in C_b(W)$. This implies that the process $(T_t,t\in [0,1])$
defined on $[0,1]\times W$ by
$$
(t,x)\to T_t(x)=x(t)+\int_0^tD_\tau\varphi(x) d\tau\,,
$$
is a Wiener process under the measure $\La d\mu$ with respect to its
natural filtration $(\calF^T_t,t\in [0,1])$, where $D_t\varphi$
represents the Lebesgue density of the map $t\to\nabla\varphi(x)(t) \in
H$ on $[0,1]$. Since $T$ is invertible, we have also 
$$
\bigvee_{t\in [0,1]}\calF_t^T=\calB(W)\,,
$$
upto $\mu$-negligeable sets. Since $\La d\mu$ is equivalent to the
Wiener measure, the process $(T_t,t\in [0,1])$ is a
$\mu$-semimartingale with respect to its natural filtration. It is
clear that it has a decomposition of the form
$$
T_t=B_t^T+A_t\,,
$$
with respect to $\mu$, where $B^T$ is a $\mu$-Brownian motion and $A$
is a process of finite variation. Since we are dealing with the
Brownian filtrations, $(A_t,t\in [0,1])$ should be absolutely
continuous with respect to the Lebesgue measure $dt$ of $[0,1]$. In
order to calculate its density it suffices to calculate the limit
$$
\lim_{h\to 0}\frac{1}{h}E\left[T_{t+h}-T_t|\calF^T_t\right]\,.
$$
To calculate this limit, it is enough to test it on the functions of
the type $g\circ T_t$:
\bea
E\left[(T_{t+h}-T_t)\,g\circ T_t\right]&=&E\left[(W_{t+h}-W_t)\,g\circ
  W_t\,L\right]\nonumber\\ 
&=&E\left[(\delta U_{[t,t+h]})g\circ W_t\,L\right]\nonumber\\
&=&E\left[(U_{[t,t+h]},\nabla(L\,g\circ W_t))_H\right]\label{l-1}\\
&=&E\left[g\circ W_t\int_t^{t+h}D_\tau Ld\tau\right]\,,\label{l-2}
\eea
where $U_{[t,t+h]}$ is the element of $H$ whose Lebesgue density is
equal to the indicator function of the interval $[t,t+h]$. Note that
for the equality (\ref{l-1}), we have used the fact that
$\delta=\nabla^\star$ under the Wiener measure $\mu$ and the equality
(\ref{l-2}) follows from the fact that the support of $\nabla(g(W_t))$
lies in the interval $[0,t]$, hence its scalar product in $H$ with
$U_{[t,t+h]}$ is zero (cf.\cite{ASU}). Hence we have 
\beaa
\lim_{h\to
  0}\frac{1}{h}E\left[T_{t+h}-T_t|\calF^T_t\right]&=&-E[D_tf\circ
T|\calF^T_t]\\
&=&-E_\nu[D_tf|\calF_t]\circ T\,,
\eeaa
$dt\times d\mu$-almost surely, where the last inequality follows from
the fact that $T^{-1}\left(\calF_t\right)= \calF_t^T$. Hence we have
proven 
\begin{proposition}
\label{semi-1}
The transport  process $(T_t,t\in [0,1])$ is a
$(\mu,(\calF_t^T))$-semimartingale with its canonical decomposition
\beaa
T_t&=&B^T_t-\int_0^t E_\nu[D_\tau f|\calF_t]\circ T\,\,d\tau\\
&=&B^T_t-\int_0^t E[D_\tau f\circ T|\calF_t^T]\,\,d\tau\,.
\eeaa
\end{proposition}
\noindent
We can give now the It\^o solution of the Monge-Amp\`ere equation:
\begin{theorem}
\label{M-A-Ito}
Assume that $f\in \DD_{2,1}$ be such that $c=E[\exp(-f)]<\infty$, denote
by $L$ the probability density defined by $\frac{1}{c}e^{-f}$ and by
$\nu$ the probability $d\nu=Ld\mu$. Assume that $d_H(\mu,\nu)<\infty$
and let $T=I_W+\nabla\varphi$ be the transport map whose properties
are announced in Theorem \ref{gaussian-case}. We have then 
\begin{equation}
\label{La-Ito}
\La=\exp\left\{\int_0^1E_\nu[D_tf|\calF_t]\circ
  TdB^T_t-\frac{1}{2}\int_0^1 E_\nu[D_tf|\calF_t]^2\circ
  T\,dt\right\}\,.
\end{equation}
\end{theorem}
\proof
From the It\^o representation formula \cite{ASU-0}, we have 
$$
L=\exp\left\{-\int_0^1E_\nu[D_tf|\calF_t]dW_t-\frac{1}{2}\int_0^1
  E_\nu[D_tf|\calF_t]^2\,dt\right\}\,. 
$$
Since the Girsanov measure for $T$ has the density $\La$ given by 
$$
\La=\frac{1}{L\circ T}\,,
$$
we have, using  the identity $T^{-1}(\calF_t)=\calF^T_t$ and
Proposition \ref{semi-1}, 
\beaa
L\circ T &=&\exp\left\{-\int_0^1E_\nu[D_tf|\calF_t]\circ
  TdT_t-\frac{1}{2}\int_0^1  E_\nu[D_tf|\calF_t]^2\circ
  T\,dt\right\}\\
&=&\exp\left\{-\int_0^1E_\nu[D_tf|\calF_t]\circ
  T\left(dB^T_t-E_\nu[D_tf|\calF_t]\circ T dt\right)\right.\\
&&\left.-\frac{1}{2}\int_0^1
  E_\nu[D_tf|\calF_t]^2\circ T\,dt\right\}\,, 
\eeaa
which is exactly the inverse of the expression given by the relation
(\ref{La-Ito}). 
\qed

\noindent
The following proposition explains the relation between   the semimartingale
representation of $T$ and the polar factorization studied in Section
\ref{pol-fac}:
\begin{proposition}
\label{pol-sm}
Let  $X$ be  the process defined by
$$
X_t=W_t+\int_0^tE_\nu[D_\tau f|\calF_\tau]d\tau\,, 
$$
then $T\circ X$ is a $\nu$-rotation, i.e., $T\circ X(\nu)=\nu$, in
fact it is the minimal $\nu$-rotation in the sense that 
$$
\inf_{\calO\in \calR_\nu}E_\nu[|O-X|_H^2]=E_\nu[|T\circ X-X|_H^2]\,,
$$
where $\calR_\nu$ denotes the set of transformations preserving the
measure $\nu$. Finally the Brownian motion $B^T$ is the rotation
corresponding to $X\circ T$.
\end{proposition}
\proof
Since $E[L]=1$, $\nu$ is the Girsanov measure for the transformation
$X$, consequently, we have 
\beaa
E_\nu[g(T\circ X)]&=&E[g(T\circ X)\,L]\\
&=&E[g(T)]\\
&=&E_\nu[g]\,,
\eeaa
for any $g\in C_b(W)$ and this implies $T\circ X(\nu)=\nu$. Let now
$\calO\in \calR_\nu$, then the measure $\calO\times X(\nu)$ belongs to
$\Sigma(\nu,\mu)$. Since $T\times I(\mu)$ is the solution of MKP in
$\Sigma(\nu,\mu)$, we have 
$$
E_\nu[|\calO-X|_H^2]\geq E_\nu[|T\circ X-X|_H^2]=d_H(\mu,\nu)^2\,.
$$
The uniqueness follows from the same argument as used in the proof of 
Theorem \ref{polar-fac-thm}. The last claim is obvious since $X\circ
T$ is a $\mu$-rotation, hence as a process it is a Brownian motion,
then by comparing it with the result of Proposition \ref{semi-1}, we
see that $B^T=X\circ T$.
\qed

\noindent
Let us give some immediate consequences of these results whose proof
follows immediately from the results of this section and from Theorem
\ref{jacobian-calcul} :
\begin{corollary}
\label{calcul-cor}
We have the following identity
\beaa
-\log E[e^{-f}]&=&E\left[f\circ T+\frac{1}{2}\int_0^1
E_\nu[D_tf|\calF_t]^2\circ T\,dt\right]\\
&=&E\left[f\circ T+\frac{1}{2}\int_0^1
E[D_tf\circ T|\calF^T_t]^2\,dt\right]\,.
\eeaa
If, furthermore, $f$ is $H$-convex, then we also have 
$$
-\log E[e^{-f}]=E\left[f\circ
  T-\log\dett(I_H+\nabla^2\varphi)+\frac{1}{2}|\nabla\varphi|_H^2\right]\,.
$$
In particular we have the exact characterization of the Wasserstein
distance between $\mu$ and $\nu$:
\beaa
\frac{1}{2}d_H^2(\mu,\nu)&=&
E\left[\log\dett(I_H+\nabla^2\varphi)\right]+\frac{1}{2}
E\left[\int_0^1E_\nu[D_tf|\calF_t]^2\circ T\,dt\right]\\ 
&=&E\left[\log\dett(I_H+\nabla^2\varphi)\right]+E[L\log L]\,.
\eeaa
\end{corollary}

\noindent{\small{
\begin{itemize}
\item D. Feyel, Universit\'e d'Evry-Val-d'Essone, 91025 Evry Cedex, France.
E-mail: feyel@maths.univ-evry.fr
\item A. S. \"Ust\"unel, ENST, D\'ept. Infres, 46, rue Barrault, 75634
  Paris Cedex 13, France.
E-mail: ustunel@enst.fr
\end{itemize}
}}

\end{document}